\documentclass[12pt]{article}
\usepackage{amsfonts,amsmath,latexsym,amssymb,mathrsfs}
\usepackage{graphicx}
\usepackage{hyperref}
\usepackage{wrapfig}
\usepackage{breakcites}
\usepackage{color}
\usepackage{breakcites}
\usepackage{geometry}
\usepackage{enumerate}
\usepackage[T1]{fontenc}
\usepackage[english]{babel}
\usepackage{dsfont,bm,bbm,mathrsfs}
\usepackage{stmaryrd}
\usepackage{url,color}
\usepackage{textcase}

\usepackage{tikz}
\usetikzlibrary{intersections}

\allowdisplaybreaks

\newtheorem{thm}{Theorem}[section]
\newtheorem{prop}[thm]{Proposition}

\newtheorem{lma}[thm]{Lemma}
\newtheorem{cor}[thm]{Corollary}
\newtheorem{exam}[thm]{Example}
\newtheorem{countexam}[thm]{Counterexample}
\newtheorem{rem}[thm]{Remark}
\newtheorem{con}[thm]{Conjecture}

\makeatletter
\def\given{\mskip 0.5mu plus 0.25mu\vert\mskip 0.5mu plus 0.15mu}
\newcounter{@bracketlevel}
\def\@bracketfactory#1#2#3#4#5#6{
\expandafter\def\csname#1\endcsname##1{%
\addtocounter{@bracketlevel}{1}%
\global\expandafter\let\csname @middummy\alph{@bracketlevel}\endcsname\given%
\global\def\given{\mskip#5\csname#4\endcsname\vert\mskip#6}\csname#4l\endcsname#2##1\csname#4r\endcsname#3%
\global\expandafter\let\expandafter\given\csname @middummy\alph{@bracketlevel}\endcsname
\addtocounter{@bracketlevel}{-1}}%
}
\def\bracketfactory#1#2#3{%
\@bracketfactory{#1}{#2}{#3}{relax}{0.5mu plus 0.25mu}{0.5mu plus 0.15mu}
\@bracketfactory{b#1}{#2}{#3}{big}{1mu plus 0.25mu minus 0.25mu}{0.6mu plus 0.15mu minus 0.15mu}
\@bracketfactory{bb#1}{#2}{#3}{Big}{2.4mu plus 0.8mu minus 0.8mu}{1.8mu plus 0.6mu minus 0.6mu}
\@bracketfactory{bbb#1}{#2}{#3}{bigg}{3.2mu plus 1mu minus 1mu}{2.4mu plus 0.75mu minus 0.75mu}
\@bracketfactory{bbbb#1}{#2}{#3}{Bigg}{4mu plus 1mu minus 1mu}{3mu plus 0.75mu minus 0.75mu}
}
\bracketfactory{clc}{\lbrace}{\rbrace}
\bracketfactory{clr}{(}{)}
\bracketfactory{cls}{[}{]}
\bracketfactory{abs}{\lvert}{\rvert}
\bracketfactory{norm}{\Vert}{\Vert}
\bracketfactory{floor}{\lfloor}{\rfloor}
\bracketfactory{ceil}{\lceil}{\rceil}
\bracketfactory{angle}{\langle}{\rangle}

\def\IP{\prob}

\makeatother

\newcommand{\law}{\mathscr{L}}

\newcommand{\equald}{\stackrel{d}{=}}
\newcommand{\Pro}{\mathbb{P}} 
\newcommand{\prob}{\Pro}
\DeclareMathOperator{\E}{\mathbb{E}}

\newcommand{\R}{\mathbb{R}}
\newcommand{\N}{\mathbb{N}}

\newcommand{\Z}{\mathbb{Z}}

\DeclareMathOperator{\Var}{\mathrm{Var}}
\DeclareMathOperator{\var}{\mathrm{Var}}

\usepackage{subfig}

\newcommand{\bone}{{\bf 1}}

\newcommand{\dk}{{d_{\rm K}}}
\newcommand{\dw}{{d_{\rm W}}}

\newcommand{\Pn}{{\rm Pn}}
\newcommand{\NB}{{\rm NB}}
\newcommand{\CP}{{\rm CP}}

\def\a{\alpha}

\def\l{\lambda}

\def\dw{d_{\mathrm{W}}}

\def\tpi{\overline{\tau^+_i}}
\def\tpi1{\overline{\tau^+_{i-1}}}

\makeatletter
 
\newcommand{\qed}{\nopagebreak\hspace*{\fill}
{\vrule width6pt height6ptdepth0pt}\par}

\def\ignore#1{}

\makeatletter
\newcommand*\labelcounter[2]{\begingroup
  \protected@edef\@currentlabel{\csname p@#1\endcsname\csname the#1\endcsname}%
  \label{#2}\endgroup}
\newcommand*\refsetcounter[2]{\setcounter{#1}{#2}%
  \protected@edef\@currentlabel{\csname p@#1\endcsname\csname the#1\endcsname}%
  }
\makeatother

\newcounter{rtaskno}

\ignore{\newcounter{con}
}

\newcounter{con}%for constants

\numberwithin{equation}{section}

\title{\sc\large%\MakeUppercase
{Geometric sums, size biasing and zero biasing
}}

\author{ Qingwei Liu\footnote{School of Mathematics and Statistics,
The University of Melbourne,
VIC 3010, Australia, E-mail: qingweil@student.unimelb.edu.au. Work supported in part by China Scholarship Council.}
\ and \
 Aihua Xia\footnote{School of Mathematics and Statistics,
The University of Melbourne,
VIC 3010, Australia, E-mail: aihuaxia@unimelb.edu.au. Work supported by Australian Research Council Grant No DP190100613.}
}

\def\parsedate #1:20#2#3#4#5#6#7#8\empty{20#2#3-#4#5-#6#7}
\def\moddate{\expandafter\parsedate\pdffilemoddate{\jobname.tex}\empty}
\date{\moddate}

\begin{document}
\maketitle

\begin{abstract}
The geometric sum plays a significant role in risk theory and reliability theory \cite{Kala97} and a prototypical example of the geometric sum is R\'enyi's theorem~\cite{Renyi56} saying a sequence of suitably  parameterised geometric sums converges to the exponential distribution.  There is extensive study of the accuracy of exponential distribution approximation to the geometric sum \cite{Sugakova95,Kala97,PekozRollin11} but there is little study on its natural counterpart of gamma distribution approximation to negative binomial sums. 
In this note, we show that a nonnegative random variable follows a gamma distribution if and only if its size biasing equals its zero biasing. We combine this characterisation with Stein's method to establish simple bounds for gamma distribution approximation to the sum of nonnegative independent random variables, a class of compound Poisson distributions and the negative binomial sum of random variables.
\end{abstract}

\section{Introduction}
R\'enyi's theorem~\cite{Renyi56} states that $\mathscr{L}(p\sum_{i=1}^NX_i)\to\rm{Exp}(1)$ as $p\to0$, where $\mathscr{L}$ denotes the distribution, $\{X_i\}$ is a sequence of independent and identically distributed (\textit{i.i.d.}) random variables with $\E X_1=1$, $N$ is a geometric random variable with distribution $\prob(N=k)=p(1-p)^k$ for $k\in\Z_+:=\{0,1,2,{\dots}\}$,
denoted by $N\sim{\rm Ge}(p)$, and $N$ is independent of $\{X_i\}$. Geometric sums are a natural object in risk theory and reliability theory \cite{Kala97} and, under mild conditions, they are asymptotically close to the exponential distribution (see \cite{Sugakova95,Kala97} and references therein). The accuracy of exponential distribution approximation can be estimated through renewal techniques hinged on the memoryless property of the geometric distribution \cite{Sugakova95,Kala97} but these techniques seem to be less efficient in tackling its natural counterpart of gamma distribution approximation to negative binomial sums. Stein's method related to gamma distribution approximation \cite{DZ91,Luk94,PekozRollin11,GPR17,Gaunt19,Slepov21} is more flexible than renewal techniques and 
gamma distribution approximation of random variables in a fixed Wiener chaos of a general Gaussian process using the Malliavin-Stein method has been investigated by \cite{ET15,LNP15,NourdinPeccati09}.
The key to the success of using Stein's method is to find suitable distributional transformations of the random object under consideration \cite{GoldsteinReinert05}, hence the first obstacle we need to overcome is to find such distributional transformations characterising the gamma distribution. To this end, let us recall the two most commonly used distributional transformations, namely, size biasing and zero biasing.

For any nonnegative random variable $V$ with finite mean $\mu$, we say that $V^s$ has $V$-size biased distribution if 
\begin{equation}\label{defsizebias}
\E[Vf(V)]=\mu\E f(V^s)
\end{equation}
for all functions $f$ such that $\E|Vf(V)|<\infty$. For the geometric sum $S:=\sum_{i=1}^NX_i$ mentioned above \cite{Renyi56}, it is a routine exercise to verify that $S^s\equald S+S'+X_1^s$,
 where $\equald$ stands for ``equal in distribution'', $S'\equald S$, and $S,S',X_1^s$ are independent. This form of size biasing does not seem promising to study gamma distribution approximation of the geometric sum using the size biasing only.

For a random variable $W$ with mean $\mu$ and variance $\sigma^2\in(0,\infty)$, we say $W^z$ has the $W$-zero biased distribution if for all differentiable function $f$ with $\E |Wf(W)|<\infty$ ,
\begin{equation}\label{defzerobias}
\E [(W-\mu)f(W)]=\sigma^2 \E f'(W^z).\end{equation}
Zero biased distribution was first introduced in \cite{GoldsteinReinert97} inspired by the following observation in \cite{Stein72}: $Z$ is a normal random variable with zero mean and variance $\sigma^2$ if and only if for all absolutely continuous $f$ with $\E|Zf(Z)|<\infty$,
$$\E[Zf(Z)]=\sigma^2\E f'(Z).$$ 
The discrete version of zero biasing was introduced in \cite{Goldstein06} and it is slightly different from the zero biasing defined above.
It is an elementary exercise to verify that zero biasing satisfies $W^z\overset{d}{=}(W-a)^z+a$ for all $a\in\R$.  

As observed in \cite{GoldsteinReinert05}, both size biasing and zero biasing are special cases of the transformation that can be used to construct approximation theory based on a distribution which is the fixed point of the transformation. Both of the biasing transformations play significant roles in a wide range of distributional approximations including normal~\cite{CGS}, Poisson~\cite{BHJ} and exponential~\cite{PekozRollin11}. 

In Section~\ref{chractgamma}, we show that the gamma distribution is uniquely characterised by the property that its size biased distribution is the same as its zero biased distribution. We then combine this characterisation with Stein's method in Section~\ref{approximationtheory} to establish simple bounds for gamma distribution approximation with application to the sum of independent nonnegative random variables. As the gamma distribution is in the family of infinitely divisible distributions, we present the direct relationship between size biasing and zero biasing for infinitely divisible distributions on $\R_+:=[0,\infty)$ in Lemma~\ref{IDzerobiasing}. In the remaining part of Section~\ref{App01}, we consider gamma distribution approximation to a class of compound Poisson distributions and the negative binomial sum of random variables. In particular, our result provides an intuitive explanation why the gamma distribution is a continuous counterpart of the negative binomial distribution and why the gamma process is a pure jump increasing process while the gamma distribution is continuous. 

\section{Gamma distribution: the intersection of size biasing and zero biasing}\label{chractgamma}
 The gamma distribution $\Gamma(r,\a)$ we consider has the probability density function $\frac1{\Gamma(r)}\a^rx^{r-1}e^{-\a x}$, $x>0$. The following theorem states that the unique distribution having the same size biased distribution and zero biased distribution is a gamma distribution. A well-known fact is that for $W\sim\Gamma(r,\a)$,  $W^s\sim\Gamma(r+1,\a)$ \cite{AGK19}.
 
 \begin{thm}\label{biasidentity}
 For a random variable $W\ge 0$ with mean $\mu$ and variance $\sigma^2\in(0,\infty)$, we have $W\sim\Gamma(r,\a)$ for some $r,\a>0$ if and only if 
 \begin{equation}\label{criteria}W^s\equald W^z.\end{equation}
 \end{thm}

\noindent{\bf Proof}~The proof relies on Stein's identity of the gamma distribution \cite{Luk94}:
\begin{equation}\label{gammaid}x f''(x)+(r-\a x)f'(x)=h(x)-\Gamma_{r,\a}h,\end{equation} 
where $\Gamma_{r,\a}h:=\E h(Z)$ for $Z\sim\Gamma(r,\a)$. Therefore, the random variable $W\sim\Gamma(r,\a)$ if and only if 
\begin{equation*}\E[Wf''(W)+(r-\a W)f'(W)]=0
\end{equation*}
for all twice differentiable functions $f$ such that the expectations $\E|Wf''(W)|$ and $\E|Wf'(W)|$ are finite.

For the necessity part, assume that $W\sim\Gamma(r,\a)$ for some $r,\a>0$, hence $$\mu=\frac r\a,~~~\sigma^2=\frac r{\a^2}.$$
For a differentiable $h$ such that $\E|Wh(W)|<\infty$, 
 \begin{eqnarray*}
\sigma^{-2}\E[(W-\mu)h(W)]&=&\a\E[h(W^s)-h(W)]\\
&=&\a\int_0^{\infty}h(x)\left[\frac1{\Gamma(r+1)}\a^{r+1}x^{r}e^{-\a x}-\frac1{\Gamma(r)}\a^rx^{r-1}e^{-\a x}\right]dx\\
&=&-\int_0^{\infty}h(x)d\left[\frac1{\Gamma(r+1)}\a^{r+1}x^{r}e^{-\a x}\right]\\
&=&\int_0^{\infty}h'(x)\frac1{\Gamma(r+1)}\a^{r+1}x^{r}e^{-\a x}dx\\
&=&\Gamma_{r+1,\a}h',
 \end{eqnarray*} 
 which ensures that $W^z\sim\Gamma(r+1,\a)$. Since $W^s\sim\Gamma(r+1,\a)$, (\ref{criteria}) follows.
 Conversely, if (\ref{criteria}) holds, by choosing 
 $$\a:=\frac{\mu}{\sigma^2},~~~r:=\frac{\mu^2}{\sigma^2},$$ we have $\mu=\a\sigma^2$. For all twice differentiable $f$ such that the following expectations exist, we have
 \begin{eqnarray*}
&&\E[Wf''(W)+(r-\a W)f'(W)]\\
 &=&\E[Wf''(W)]-\a\E[(W-\mu)f'(W)]\\
 &=&\mu\E f''(W^s)-\a\sigma^2\E f''(W^z)\\
 &=&0.
 \end{eqnarray*}
 This ensures that $W\sim\Gamma(r,\a)$ and the proof is complete.
 
\qed

\section{Main results}\label{approximationtheory}

In this section, we bound the errors of gamma distribution approximation in terms of the Wasserstein distance $\dw$ and the Kolmogorov distance $\dk$ defined as
\begin{eqnarray}
\label{wass01}\dw(\law(X),\law(Y))&:=&\sup_{f\in\mathcal{F}_W}|\E f(X)-\E f(Y)|,\\
\dk(\law(X),\law(Y))&:=&\sup_{f\in\mathcal{F}_K}|\E f(X)-\E f(Y)|,\nonumber
\end{eqnarray}
where $\mathcal{F}_W:=\{f:\R\to\R,|f(x)-f(y)|\le|x-y|,\forall x,y\in\R\}$ and $\mathcal{F}_K:=\{\bone_{\{\cdot\le z\}},\forall z\in\R\}$.

 \begin{thm}\label{thm1}
 Let $W$ be a nonnegative random variable with mean $\mu$ and variance $\sigma^2\in(0,\infty)$. For $\a=\mu/\sigma^2$, $r=\mu^2/\sigma^2$,
\begin{align}
\dw(\law(W),\Gamma(r,\a))&\le8\sqrt{\frac{3\mu}{r+2}}\Theta^{1/2}+\frac{8r}{r+2}\Theta,\label{thm1.1}\\
 \dk(\law(W),\Gamma(r,\a))&\le
 a_{r,\a}\Theta^{\frac{r\wedge 1}{r\wedge 1+2}}+b_{r,\a}\Theta^{\frac{r\wedge 1+1}{r\wedge 1+2}},\label{thm1.2}
 \end{align}
 where $\Theta:=\dw(\law(W^s),\law(W^z))$ and
 \begin{align}
 a_{r,\a}&=\left\{\begin{array}{ll}
 (1/2+1/r)\left(\frac{48\mu}{r+2}\right)^{\frac{r}{r+2}}\left[\frac{\a^r}{\Gamma(r)}\right]^{\frac2{r+2}},&0<r<1,\\
 3\left(\frac{6\mu}{r+2}\right)^{1/3}\left[\frac{\a}{\Gamma(r)}\left(\frac{r-1}e\right)^{r-1}\right]^{2/3},&r\ge 1;\end{array}\right.\label{thm1.2-1}\\
 b_{r,\a}&=\left\{\begin{array}{ll}8\a\left(\frac{\mu}{r+2}\right)^{\frac{r+1}{r+2}}\left[\frac{\a^r}{48\Gamma(r)}\right]^{\frac1{r+2}},&0<r<1,\\
 4\a\left(\frac{\mu}{r+2}\right)^{2/3}\left[\frac{\a}{6\Gamma(r)}\left(\frac{r-1}e\right)^{r-1}\right]^{1/3},&r\ge 1.\end{array}\right.
 \label{thm1.2-2}
 \end{align}
 \end{thm}
 
\noindent{\bf Proof}
For any $h\in\mathcal{F}_W$ and $\delta>0$, we can construct a smooth interpolating spline function $\tilde{h}$ as follows. For $i\in\mathbb{Z}:=\{0,\pm1,\pm2,\cdots\}$, let $x_i=-i\delta$, $(x_i,y_i):=(x_i,h(x_i))$ and 
$$
 \tilde{h}(x):=y_i+(y_{i+1}-y_i)\phi\left(\frac{x-x_i}{x_{i+1}-x_i}\right),\quad x\in[x_i,x_{i+1}),
$$
where 
$$
\phi(t)=\left\{
\begin{aligned}
&0, \quad &t<0,\\
&2t^2, \quad &0\le t<1/2,\\
&1-2(1-t)^2, \quad &1/2\le t\le1,\\
&1, \quad &t>1.
\end{aligned}\right.
$$
It is easy to verify that $\tilde{h}$ is smooth and 
$$
\|\tilde{h}'\|:=\sup_{x\in \R}|\tilde{h}'(x)|\le2,\quad \|h-\tilde{h}\|\le\delta,\quad \|\tilde{h}''\|\le4\delta^{-1},
$$
since $\left|\frac{y_{i+1}-y_i}{x_{i+1}-x_i}\right|\le1$.
Let $\tilde{f}:=f_{\tilde{h}}$ be the solution of the Stein equation (\ref{gammaid}) with $\tilde{h}$ in place of $h$, then
\begin{equation}\|\tilde{f}'''\|\le \frac{2}{r+2}(3\|\tilde{h}''\|+2\a\|\tilde{h}'\|)\le\frac{8}{r+2}\left(3\delta^{-1}+\a\right),\label{Thm1.1proof1}
\end{equation}
see \cite[Theorem~2.1]{GPR17}. By virtue of (\ref{gammaid}), (\ref{defsizebias}) and (\ref{defzerobias}), we have  
\begin{eqnarray*}
\E \tilde{h}(W)-\Gamma_{r,\a}\tilde{h}&=&\E[W\tilde{f}''(W)+(r-\a W)\tilde{f}'(W)]\\
&=&\mu\E \tilde{f}''(W^s)-\a\E[(W-\mu)\tilde{f}'(W)]\\
&=&\mu\E[\tilde{f}''(W^s)-\tilde{f}''(W^z)],
\end{eqnarray*}
which, together with (\ref{Thm1.1proof1}), implies
\begin{equation}
|\E \tilde{h}(W)-\Gamma_{r,\a}\tilde{h}|\le\frac{8\mu}{r+2}\left(3\delta^{-1}+\a\right)\Theta.\label{thmproof1}
\end{equation}
By the triangle inequality, we have
 \begin{eqnarray*}
 |\E h(W)-\Gamma_{r,\a}h|&\le&|\E h(W)-\E\tilde{h}(W)|+|\E\tilde{h}(W)-\Gamma_{r,\a}\tilde{h}|+|\Gamma_{r,\a}\tilde{h}-\Gamma_{r,\a}h|\\
 &\le&2\delta+\frac{8\mu}{r+2}\left(3\delta^{-1}+\a\right)\Theta,
 \end{eqnarray*}
and (\ref{thm1.1}) follows from (\ref{wass01}) and $\delta=\sqrt{\frac{12\mu}{r+2}}\Theta^{1/2}$. 

The proof of (\ref{thm1.2}) relies on the following concentration inequality of $Z\sim\Gamma(r,\a)$: for $\delta>0$ and $z\ge0$,
$$
\IP(z<Z\le z+\delta)\le\epsilon(\delta):=\left\{
\begin{aligned}
&\frac{\a^r\delta^r}{\Gamma(r+1)},\quad &0<r<1,\\
&M(r,\a)\delta,\quad &r\ge1,
\end{aligned}\right.
$$
where $M(r,\a):=\frac{\a}{\Gamma(r)}\left(\frac{r-1}{e}\right)^{r-1}$ is the maximum of the density function of $\Gamma(r,\a)$.
In fact, for $r\ge1$, the bound is obvious, and for $r\in(0,1)$, the bound follows from the fact that
$g(\delta):=\int_z^{z+\delta}x^{r-1}e^{-\a x}dx-\delta^r/r$
 is decreasing in $\delta$ and $g(0)=0$. 

Assume now $h_z(\cdot)=\bone_{\{\cdot\le z\}}\in\mathcal{F}_K$, $z\ge0$. Denote 
$$
\tilde{h}_z(x):=1-\phi\left(\frac{x-z}{\delta}\right),
$$ and note that $\|\tilde{h}_z'\|\le2\delta^{-1}$ and $\|\tilde{h}_z''\|\le4\delta^{-2}$. As 
\begin{eqnarray*}
\IP(W\le z)-\IP(Z\le z)&\le&\E\tilde{h}_z(W)-\E\tilde{h}_z(Z)+\E\tilde{h}_z(Z)-\IP(Z\le z)\\
&\le&\E\tilde{h}_z(W)-\E\tilde{h}_z(Z)+\IP(z<Z\le z+\delta)\\
&\le&|\E \tilde{h}_z(W)-\Gamma_{r,\a}\tilde{h}_z|+\epsilon(\delta),
\end{eqnarray*} and, following the same argument for (\ref{thmproof1}), we have
$$
|\E \tilde{h}_z(W)-\Gamma_{r,\a}\tilde{h}_z|\le\mu\|\tilde{f}_z'''\|\Theta\\
\le\frac{8\mu}{r+2}(3\delta^{-2}+\a\delta^{-1})\Theta=:e(\delta).
$$
Analogously,
\begin{eqnarray*}
\IP(W\le z)-\IP(Z\le z)&\ge&\E\tilde{h}_{z-\delta}(W)-\E\tilde{h}_{z-\delta}(Z)-\IP(z-\delta<Z\le z)\\
&\ge&-e(\delta)-\epsilon(\delta).
\end{eqnarray*}
Therefore, for $0<r<1$, we have
$$
\dk(\law(W),\Gamma(r,\a))\le e(\delta)+\frac{\a^r\delta^r}{\Gamma(r+1)}
=\frac{8\mu}{r+2}(3\delta^{-2}+\a\delta^{-1})\Theta+\frac{\a^r\delta^r}{\Gamma(r+1)},
$$
 which implies (\ref{thm1.2}) with $\delta=\left[\frac{48\mu \Theta\Gamma(r)}{(r+2)\a^r}\right]^{1/(r+2)}$.
Similarly, for $r\ge1$, 
\begin{eqnarray*}
\dk(\law(W),\Gamma(r,\a))\le\frac{8\mu}{r+2}(3\delta^{-2}+\a\delta^{-1})\Theta+M(r,\a)\delta,
\end{eqnarray*}
which, together with $\delta=\left(\frac{48\mu \Theta}{(r+2)M(r,\a)}\right)^{1/3}$, ensures (\ref{thm1.2}).
\qed
 
 Theorem~\ref{thm1} says that the accuracy of gamma distribution approximation with respect to $\dw$ and $\dk$ is determined by $\Theta=\dw(\law(W^s),\law(W^z))$. The next corollary says that for the sum of nonnegative independent random variables, $\Theta$ can be easily bounded.
 
 \begin{cor}\label{indeptsum} 
 Let $\{X_i:\ 1\le i\le n\}$ be nonnegative independent random variables with positive finite variances and $W=\sum_{i=1}^nX_i$,
 then (\ref{thm1.1}) and (\ref{thm1.2}) hold with
 $$\Theta=\E|X_{I_1}^s+X_{I_2}-X_{I_1}-X_{I_2}^z|,$$
where $X_i^s$ and $X_i^z$ have the size-biased distribution and zero-biased distribution of $X_i$, respectively, $I_1,I_2$ are random indices, independent of $X_1,\dots,X_n$, with distributions
\begin{equation}\IP(I_1=i)=\frac{\E X_i}{\E W},~~~\IP(I_2=i)=\frac{\Var(X_i)}{\Var(W)},\label{index1}\end{equation}
for $i=1,\dots,n$.
In particular, when $X_1,\dots,X_n$ are i.i.d. random variables, we have  (\ref{thm1.1}) and (\ref{thm1.2}) with
$ \Theta= \E|X_1^s-X_1^z|.
 $
 \end{cor}

\noindent{\bf Proof} According to \cite{GoldsteinRinott05} and \cite{GoldsteinReinert97}, we can set 
$W^s:=\sum_{j\ne I_1}X_j+X_{I_1}^s$ and $W^z:=\sum_{j\ne I_2}X_j+X_{I_2}^z$, hence 
\begin{equation}
\E|W^s-W^z|%&=&\E|(W-X_{I_1}+X_{I_1}^s)-(W-X_{I_2}+X_{I_2}^z)|\non\\
=\E|X_{I_1}^s+X_{I_2}-X_{I_1}-X_{I_2}^z|.\label{noniid}
\end{equation}
When $X_1,\dots,X_n$ are i.i.d. random variables, we may take $I_1=I_2$, which is uniformly distributed on $\{1,\dots,n\}$. Hence, (\ref{noniid}) can be reduced to
\begin{equation*}\E|X_{I_1}^s-X_{I_1}^z|=\E|X_1^s-X_1^z|,\end{equation*}
as claimed.
\qed

Convolutions of gamma distributions commonly arise in statistics \cite{VU09,CovoElalouf14}. The next corollary quantifies gamma distribution approximation to such convolutions.
 
\begin{cor}%\label{indeptgamma}
If $X_i\sim\Gamma(r_i,\a_i)$, $1\le i\le n$, are independent gamma distributed random variables, with $\mu=\sum_{i=1}^nr_i/\a_i$, $\sigma^2=\sum_{i=1}^nr_i/\a_i^2,$ $\a=\mu/\sigma^2$, $r=\a\mu$, then  (\ref{thm1.1}) and (\ref{thm1.2}) hold with 
$\Theta= \E|Y_{I_1}-Y_{I_2}|$,
where $I_1$ and $I_2$ are random indices satisfying
\begin{equation}\IP(I_1=i)=\frac{r_i\a_i^{-1}}{\sum_{j=1}^nr_j\a_j^{-1}},~~~\IP(I_2=i)=\frac{r_i\a_i^{-2}}{\sum_{j=1}^nr_j\a_j^{-2}},\label{index2}\end{equation}
$Y_i\sim\Gamma(1,\a_i)$, $i=1,\dots,n,$ $\{Y_i\}$ are independent and are independent of $\{I_1,I_2\}$.
\end{cor}

\noindent{\bf Proof}~~By Corollary~\ref{indeptsum}, it suffices to bound $\E|X_{I_1}^s+X_{I_2}-X_{I_1}-X_{I_2}^z|$.
Recalling that $X_i^s\equald X_i^z\sim\Gamma(r_i+1,\a_i)$, we can set
$X_i^s:=X_i+Y_i,~~~X_i^z=X_i^s,$
where $Y_i\sim\Gamma(1,\a_i)$ is independent of $\{X_j:1\le j\le n\}$, $1\le i\le n$. The distributions of $I_1,I_2$ in (\ref{index1}) are reduced to (\ref{index2}) and $\E|X_{I_1}^s+X_{I_2}-X_{I_1}-X_{I_2}^z|= \E|Y_{I_1}-Y_{I_2}|$.
\qed

\section{Applications}\label{App01}

Before we consider applications, it is handy to have the following lemma bounding the Wasserstein distance between two gamma distributions.

\begin{lma}\label{discrepancygamma}
For $r_1,r_2,\a_1,\a_2>0$, 
\begin{equation*}\dw(\Gamma(r_1,\a_1),\Gamma(r_2,\a_2))\le\frac{|r_1-r_2|}{\a_1\vee\a_2}+(r_1\vee r_2)\left|\frac1\a_1-\frac1\a_2\right|.\end{equation*}
\end{lma}
\noindent{\bf Proof}~~In fact, $\Gamma(r,\a)$ is stochastically increasing in $r$ and stochastically decreasing in $\a$, by the triangle inequality, suppose $\a_1<\a_2$
\begin{eqnarray*}
\dw(\Gamma(r_1,\a_1),\Gamma(r_2,\a_2))&\le&\dw(\Gamma(r_1,\a_1),\Gamma(r_1,\a_2))+\dw(\Gamma(r_1,\a_2),\Gamma(r_2,\a_2))\\
&=&r_1\left|\frac1{\a_1}-\frac1{\a_2}\right|+\frac{|r_1-r_2|}{\a_2}.
\end{eqnarray*}
\qed

We will also need size biasing and zero biasing of infinitely divisible distributions on $\R_+$, see \cite[Theorem~11.2]{AGK19} and \cite[Proposition~3.8]{ArrasHoudre19} for more details. The direct relationship between the two biasings seems to be not noted anywhere in the literature so we give a proof for the relationship.

\begin{lma}\label{IDzerobiasing} Suppose $X$ is a nonnegative random variable with $\var(X)>0$.
\begin{description}
\item{(a)} $\law(X)$ is infinitely divisible if and only if there exists a random variable $X'\ge 0$ a.s. independent of $X$ such that 
\begin{equation}X^s\overset{d}{=}X+X'.\label{IDsizebiasing}\end{equation}
\item{(b)} $\law(X)$ is infinitely divisible if and only if there exists a random variable $X''$ independent of $X$ such that 
\begin{equation}\label{zeroinfi}X^z\equald X+X''.
\end{equation}
The distribution of $X''$ has the density function
\begin{equation}\label{add02}\frac1{\E X'}\IP(X'\ge x),~~~x\ge0,
\end{equation}
where $X'$ is uniquely determined in (\ref{IDsizebiasing}).
\end{description}
\end{lma}

\noindent{\bf Proof of (\ref{add02})}~The proof is a simple application of the Laplace transform. For $\theta>0$, denote $\phi_V(\theta):=\E e^{-\theta V}$ for some non-negative random variable $V$. For simplicity, we denote $\mu:=\E X$, and $\sigma^2:=\var(X)$. By taking $V=X$ and $f(x)=e^{-\theta x}$ in (\ref{defzerobias}), together with $\E X'=\E X^s-\E X=\frac{\sigma^2}{\mu}$, we have 
\begin{eqnarray*}
\phi_{X^z}(\theta)&=&-\frac1{\sigma^2\theta}\E\left[(X-\E X)e^{-\theta X}\right]\\
&=&-\frac{\mu}{\sigma^2\theta}\left[\phi_{X^s}(\theta)-\phi_X(\theta)\right]\\
&=&\phi_X(\theta)\frac{\mu}{\sigma^2\theta}\E[1-e^{-\theta X'}]\\
&=&\phi_X(\theta)\frac{\mu}{\sigma^2}\int^\infty_0e^{-\theta x}[1-\IP(X'\le x)]dx\\
&=&\phi_X(\theta)\int^\infty_0e^{-\theta x}\frac{1-\IP(X'\le x)}{\E X'}dx,
\end{eqnarray*}where the third equality is due to (\ref{IDsizebiasing}), and the fourth one is from the integration by parts. This is equivalent to (\ref{zeroinfi}).
\qed

We note that the distribution of $X''$ is also called the \textit{equilibrium distribution} with respect to $X'$ in \cite{PekozRollin11}. \cite{ArrasHoudre19} state that the only probability measure that has an additive exponential size biased distribution is the gamma distribution and \cite{PekozRollin11} say that the exponential distribution is the unique fixed point under the equilibrium transformation, their observations confirm Theorem~\ref{biasidentity} in the case of infinitely divisible distributions on $\R_+$.

The first application we consider is to estimate the difference between a compound Poisson distribution and a gamma distribution having the same mean and variance. Recall that a compound Poisson distribution, denoted by $\CP(\l,\law(X))$, is the distribution of $W=\sum_{i=1}^NX_i$, where $\{X, X_i,~i\ge1\}$ are i.i.d. random variables independent of $N\sim\Pn(\l)$. 

\begin{prop}\label{approcp}  Assume that $X\ge 0$ and $\var(X)\in(0,\infty)$ such that both $X^s$ and $X^z$ exist. 
Let $W\sim\CP(\l,\law(X))$ with $\mu=\E W$ and $\sigma^2=\var(W)$. Taking $\a=\mu/\sigma^2$, $r=\mu^2/\sigma^2$, 
then  (\ref{thm1.1}) and (\ref{thm1.2}) hold with $\Theta=\dw(\law(X^s),\law(\tilde{X}))$, where the density function of $\tilde{X}$ is given by
\begin{equation*}f_{\tilde{X}}(y)=\frac1{\E (X^2)}\E\left[ X\bone_{\{X\ge y\}}\right],~~~y\ge0.\end{equation*}
\end{prop}

\noindent{\bf Proof}~In this case, $\law(W)$ is finitely divisible and $W^s=W+X^s$ \cite[p.~7]{AGK19}, Lemma~\ref{IDzerobiasing}
ensures $W^z=W+\tilde{X}$, the claim is an immediate consequence of Theorem~\ref{thm1}. 
\qed

The intriguing phenomenon that the gamma process is a pure-jump increasing process while the gamma distribution is continuous can be well explained by the following bound. We use this example to show that the bounds in Proposition~\ref{approcp} have some room for improvement. Recall that the L\'{e}vy measure of the gamma process is $\gamma(dy)=e^{-y}/y~dy$. 

\begin{exam}\label{cpexample1}
Let $W\sim\CP(\l,\law(X))$ with $\l=\int_{\delta}^\infty\frac{e^{-x}}{x}dx$ and the density of $\law(X)$ given by $\frac1{\l}\frac{e^{-x}}{x}\bone_{\{x\ge\delta\}}$ for some $\delta>0$, then 
$$\dw(\law(W),\Gamma(1,1))\le8\sqrt{\delta}+\frac{17}3\delta.$$
\end{exam}%\Comment{need $r$ and $\alpha$}

\begin{rem}\label{Nathan}The bound in Example~\ref{cpexample1} is not of optimal order. In fact, let $\Xi$ be a Poisson point process on $(0,\infty)$ with intensity measure $\mu(dx)=\frac{e^{-x}}x~dx$, and $\tilde{\Xi}:=\Xi|_{(\delta,\infty)}$ be the restriction of $\Xi$ to $(\delta,\infty)$. Let $Z:=\sum_{\xi\in\Xi}\xi$ and $\tilde{W}:=\sum_{\xi\in\tilde{\Xi}}\xi$. It is easy to verify that $\tilde W\overset{d}{=}W$ and $Z\sim\Gamma(1,1)$, which ensure that $\law(Z)$ is stochastically bigger than $\law(W)$, hence
$$\dw(\law(W),\Gamma(1,1))=\E Z-\E W=1-e^{-\delta}\le\delta.$$
\end{rem}

\noindent{\bf Proof Example~\ref{cpexample1}}~
It is easy to see that  $$\E X=e^{-\delta}/\l,~~\E (X^2)=(1+\delta)e^{-\delta}/\l.$$
To further compute $\dw(\law(X^s),\law(\tilde{X}))$ with $\tilde{X}$ defined in Proposition~\ref{approcp}, we use the fact \cite{Val72} that for random variables $U_1$ and $U_2$ on $\R_+$,
\begin{equation*}\dw(\law(U_1),\law(U_2))=\int_0^\infty\left|\IP(U_1\ge x)-\IP(U_2\ge x)\right|dx.\end{equation*}
Therefore, it follows that 
\begin{eqnarray}
\dw(\law(X^s),\law(\tilde{X}))&=&\int_0^\infty\left\vert\IP(X^s\ge x)-\IP(\tilde{X}\ge x)\right\vert dx\nonumber\\
&=&\int_0^\infty\left\vert\frac{\E[X\bone_{\{X\ge x\}}]}{\E X}-\int_x^\infty\frac{\E[X\bone_{\{X\ge t\}}]}{\E (X^2)}dt\right\vert dx\nonumber\\
&=&\int_0^\infty\left\vert\frac{\E[X\bone_{\{X\ge x\}}]}{\E X}-\frac{\E[X(X-x)\bone_{\{X\ge x\}}]}{\E (X^2)}\right\vert dx.\label{comexam01}
\end{eqnarray}
For $0<x<\delta$, 
\begin{equation*}
\frac{\E[X(X-x)\bone_{\{X\ge x\}}]}{\E (X^2)}=\frac{\E[X(X-x)]}{\E (X^2)}=1-\frac x{1+\delta},
\end{equation*}
and for $x\ge\delta$,
\begin{eqnarray*}
\frac{\E[X(X-x)\bone_{\{X\ge x\}}]}{\E (X^2)}&=&\frac1{\E (X^2)}\int_x^\infty t(t-x)\frac1\l\frac{e^{-t}}tdt\\
&=&(1+\delta)^{-1}e^{-(x-\delta)},
\end{eqnarray*}
and 
$$\frac{\E[X\bone_{\{X\ge x\}}]}{\E X}=\left\{\begin{array}{ll}e^{-(x-\delta)},&\mbox{ for }x\ge \delta,\\
1,&\mbox{ for }0<x<\delta.\end{array}\right.$$
Combining all components together, we have from (\ref{comexam01}) that
\begin{eqnarray*}
\dw(\law(X^s),\law(\tilde{X}))&=&\int_0^\delta \frac x{1+\delta}dx+\int_\delta^\infty\frac\delta{1+\delta}e^{-(x-\delta)}dx\\
&=&\frac{\delta(1+\delta/2)}{1+\delta}\le\delta.
\end{eqnarray*}
Taking $r=\frac{e^{-\delta}}{1+\delta},\a=\frac1{1+\delta}$, by Lemma~\ref{discrepancygamma}, 
\begin{eqnarray*}
\dw(\Gamma(r,\a),\Gamma(1,1))&\le&1-\frac{e^{-\delta}}{\delta+1}+\delta\le3\delta.
\end{eqnarray*}
Hence, by the triangle inequality, we have
\begin{eqnarray*}
\dw(\law(W),\Gamma(1,1))&\le&\dw(\law(W),\Gamma(r,\a))+\dw(\Gamma(r,\a),\Gamma(1,1))\\
&\le&8\sqrt{\frac{3e^{-\delta}(1+\delta)}{e^{-\delta}+2(1+\delta)}}\sqrt{\delta}+\frac{8e^{-\delta}}{e^{-\delta}+2(1+\delta)}\delta+3\delta.
\end{eqnarray*}
Since both $\frac{e^{-\delta}(1+\delta)}{e^{-\delta}+2(1+\delta)}$ and $\frac{e^{-\delta}}{e^{-\delta}+2(1+\delta)}$ are decreasing functions of $\delta$, the claim follows. \qed 

It is tempting to ask whether for $W\sim\CP(\l,\law(X))$ with $\l\to\infty, \E W\to\frac ab$ and $\var(W)\to\frac a{b^2}$ for some $a,b>0$ are sufficient to ensure that $W$ converges in distribution to $\Gamma(a,b)$. The following example gives a negative answer to this question, indicating that $\dw(\law(X^s),\law(\tilde{X}))\to 0$ in Proposition~\ref{approcp} is also necessary for the compound Poisson distribution to be close to the gamma distribution.

\begin{countexam}%\label{counterexample1}
Let $\IP(X=1)=1/\l=1-\IP(X=0)$, then $W\sim\CP(\l,\law(X))=\Pn(1)$, $\E W=\var(W)=1$, but $W$ does not converge to $\Gamma(1,1)$ as $\lambda\to\infty$.
\end{countexam}

For $\kappa>0$ and $0<p<1$, we write $V\sim\NB(\kappa,p)$ if 
$$\IP(V=i)=\frac{\Gamma(\kappa+i)}{\Gamma(\kappa)i!}p^\kappa(1-p)^i,~~~i=0,1,\dots.$$Hence, $\E V=\kappa(1-p)/p$ and $\var(V)=\kappa(1-p)/p^2$.

To estimate gamma distribution approximation to negative binomial sums, we first bound the difference between a rescaled negative binomial distribution and a gamma distribution. 

\begin{prop}\label{cpexample2}
Let $T_p\sim\NB(\kappa,p)$, $W_p:=pT_p$. Then
\begin{align*}
\dw(\law(W_p),\Gamma(\kappa(1-p),1))&\le4\sqrt{\frac{6\kappa (1-p)p}{\kappa(1-p)+2}}+\frac{4\kappa (1-p)p}{\kappa(1-p)+2},\\\dk(\law(W_p),\Gamma(\kappa(1-p),1))&\le a_{r,\a}(0.5p)^{\frac{r\wedge 1}{r\wedge 1+2}}+b_{r,\a}(0.5p)^{\frac{r\wedge 1+1}{r\wedge 1+2}},
 \end{align*}
where  $a_{r,\a}$ and $b_{r,\a}$ are given in (\ref{thm1.2-1}) and (\ref{thm1.2-2}).
\end{prop}

\noindent{\bf Proof}
We write $T_p:=\sum_{i=1}^NX_i$, where $\{X,X_1,X_2,\dots\}$ are i.i.d. random variables with the logarithmic distribution
\begin{equation}\IP(X=i)=-\frac1{\ln (p)}\frac{(1-p)^i}i,~~~i\ge1,\label{logdist}\end{equation}
and $N\sim\Pn(-\kappa\ln (p))$, independent of $X_i$'s. It is obvious to see that $\E X=-\frac{1-p}{p\ln (p)}$, and $\E (X^2)=-\frac{1-p}{p^2\ln(p)}$, giving $$\E W_p=\var(W_p)=\kappa(1-p).$$ 
From Lemma~\ref{IDzerobiasing}, we know that $W_p=p\sum_{i=1}^NX_i$ has the size biased distribution
$$W_p^s\overset{d}{=}W_p+\left(pX\right)^s,$$where $W_p$ and $\left(pX\right)^s$ are independent. For any $i\in\N$, $x=pi$, we have
\begin{equation}\IP\left(\left(pX\right)^s=x\right)=\frac{x\IP(p X=x)}{\E(pX)}=\frac{i\IP(X=i)}{\E X}=\IP(X^s=i).\label{nbcpsb1}\end{equation}
From (\ref{logdist}), we can derive the size biased distribution of $X$: for $i\ge1$
\begin{equation}\IP(X^s=i)=\frac{i\IP(X=i)}{\E X}=p(1-p)^{i-1}.\label{nbcpsb2}
\end{equation}
Likewise, $W_p^z\overset{d}{=}W_p+\tilde{X}$, where $\tilde{X}$ is independent of $W_p$, having density
\begin{equation*}f_{\tilde{X}}(y)=\frac{\E\left[(pX)\bone_{\{pX\ge y\}}\right]}{\E \left[(pX)^2\right]}=\frac{\E\left[X\bone_{\{X\ge y/p\}}\right]}{p\E (X^2)}.\end{equation*}
Noting that for $p(i-1)<y\le p i$, we have $i-1<y/p\le i$, thus
\begin{eqnarray*}
f_{\tilde{X}}(y)&=&\frac{\E X\IP(X^s\ge y/p)}{p\E (X^2)}
=\IP(X^s\ge i)\\
&=&(1-p)^{i-1}.
\end{eqnarray*}
\begin{eqnarray*}
\dw\left(\law\left(\left(pX\right)^s\right),\law(\tilde{X})\right)&=&\int_0^\infty\left|\IP\left(\left(pX\right)^s\ge x\right)-\IP(\tilde{X}\ge x)\right|dx\\
&=&\sum_{i=1}^\infty\int_{p({i-1})}^{pi}\left|\IP\left(\left(pX\right)^s\ge x\right)-\IP(\tilde{X}\ge x)\right|dx.
\end{eqnarray*}
For $p(i-1)<x\le pi$, from (\ref{nbcpsb1}) and (\ref{nbcpsb2}), we have
\begin{equation*}
\IP\left(\left(pX\right)^s\ge x\right)=\sum_{k=i}^\infty\IP(X^s=k)=(1-p)^{i-1}.
\end{equation*}
On the other hand,
\begin{eqnarray*}
\IP(\tilde{X}\ge x)&=&\int_x^\infty f_{\tilde{X}}(y)dy
=(pi-x)(1-p)^{i-1}+\sum_{k=i}^\infty p(1-p)^k\\
&=&(pi-x)(1-p)^{i-1}+(1-p)^i
<(1-p)^{i-1}=\IP\left(\left(pX\right)^s\ge x\right).
\end{eqnarray*}

Therefore, 
\begin{eqnarray*}
\dw\left(\law\left(\left(pX\right)^s\right),\law(\tilde{X})\right)&=&\int_0^\infty\left(\IP\left(\left(pX\right)^s\ge x\right)-\IP(\tilde{X}\ge x)\right)dx\\
&=&\E\left(pX\right)^s-\E \tilde{X}\\
&=&\frac{\E\left(pX\right)^2}{\E \left(pX\right)}-\sum_{i=1}^\infty\int_{p(i-1)}^{pi}\IP(\tilde{X}\ge x)dx\\
&=&1-\sum_{i=1}^\infty\left(p(1-p)^i+0.5p^2(1-p)^{i-1}\right)\\
&=&0.5p.
\end{eqnarray*}
This, together with Proposition~\ref{approcp}, completes the proof.
\qed

\begin{cor} Let $T_p\sim\NB(\kappa,p)$ and $\{X_i\}$ be a sequence of random variables, define $S_n=\sum_{i=1}^nX_i$. Assume that $\E(X_i|T_p)=1$ and $\var(S_i|T_p=i)=i\nu^2$ for all $i\ge 1$, then $W_p=p\sum_{i=1}^{T_p}X_i$ satisfies
\begin{equation*}
\dw(\law(W_p),\Gamma(\kappa,1))\le4\sqrt{\frac{6\kappa (1-p)p}{\kappa(1-p)+2}}+\frac{4\kappa (1-p)p}{\kappa(1-p)+2}+\nu\sqrt{\kappa(1-p)p}+p\kappa.
\end{equation*}
\end{cor}

\noindent{\bf Proof}~~By the triangle inequality,
\begin{align*}
&\dw(\law(W_p),\Gamma(\kappa,1))\\
&\le\dw(\law(W_p),\law(pT_p))+\dw(\law(pT_p),\Gamma(\kappa(1-p),1))+\dw(\Gamma(\kappa(1-p),1),\Gamma(\kappa,1))\\
%&\le\E|W_p-pT_p|+\dw(\law(pT_p),\Gamma(\kappa(1-p),1))+\dw(\Gamma(\kappa(1-p),1),\Gamma(\kappa,1))\\
&\le\E|W_p-pT_p|+4\sqrt{\frac{6\kappa p(1-p)}{\kappa(1-p)+2}}+\frac{4\kappa p(1-p)}{\kappa(1-p)+2}+p\kappa,
\end{align*}
where the last inequality comes from Proposition~\ref{cpexample2} and Lemma~\ref{discrepancygamma} using the fact that $\Gamma(r,1)$ is stochastically increasing in $r$.
The remaining part $\E|W_p-pT_p|$ is bounded by the Cauchy-Schwarz inequality:
\begin{align*}
\E|W_p-pT_p|&=p\E\left|\sum_{i=1}^{T_p}(X_i-1)\right|
%&\le&p\|h'\|\sqrt{\E\left(\sum_{i=1}^{T_p}(X_i-1)\right)^2}\\
\le p\sqrt{\var\left(\sum_{i=1}^{T_p}(X_i-1)\right)}\\
&=p\sqrt{\E\left[\var\left(\sum_{i=1}^{T_p}(X_i-1)\bigg\vert T_p\right)\right]}
=p\sqrt{\E[T_p\nu^2]}
=\nu\sqrt{\kappa(1-p)p}.
\end{align*}
\qed
\noindent\textbf{Acknowledgements} We thank Nathan Ross for suggesting the direct coupling proof in Remark~\ref{Nathan} and in the context of infinitely divisible distributions on $\R_+$, Theorem~\ref{biasidentity} is indirectly confirmed in the literature.

\nocite{*}

%%%%%%%%%%%%%%%%%%%%%%%%%%%%%
%%%%%%% References   %%%%%%%%
%%%%%%%%%%%%%%%%%%%%%%%%%%%%%

\def\ac{{Academic Press}~}
\def\aap{{Adv. Appl. Prob.}~}
\def\ap{{Ann. Probab.}~}
\def\anap{{Ann. Appl. Probab.}~}
\def\eljp{{Electron. J. Probab.}~} 
\def\eljs{{Electron. J. Stat.}~}
\def\jap{{J. Appl. Probab.}~}
\def\jws{{John Wiley~$\&$ Sons}~}
\def\ny{{New York}~}
\def\ptrf{{Probab. Theory Related Fields}~}
\def\sp{{Springer}~}
\def\spa{{Stochastic Process. Appl.}~}
\def\sv{{Springer-Verlag}~}
\def\tpa{{Theory Probab. Appl.}~}
\def\zw{{Z. Wahrsch. Verw. Gebiete}~}


\begin{thebibliography}{9}
\bibitem[Arras \& Houdr\'e~(2019)]{ArrasHoudre19} Arras, B. \& Houdr\'e, C.~(2019). On Stein's Method for Infinitely Divisible Laws with Finite First Moment. \emph{Springer International Publishing}.

\bibitem[Arratia, Goldstein \& Kochman~(2019)]{AGK19} Arratia, R., Goldstein, L. \& Kochman, F.~(2019). Size bias for one and all. \emph{Probab. Surveys} \textbf{16}, 1--61. \href{https://mathscinet.ams.org/mathscinet-getitem?mr=MR3896143}{MR3896143}

\bibitem[Barbour, Holst \& Janson~(1992)]{BHJ} Barbour, A. D., Holst, L. \& Janson, S.~(1992). Poisson Approximation. \emph{The Clarendon Press Oxford University Press}.

\bibitem[Chen, Goldstein \& Shao~(2010)]{CGS} Chen, L. H. Y., Goldstein, L. \& Shao, Q. M.~(2010). Normal approximation by Stein's method. \emph{Springer Science \& Business Media.}

\bibitem[Covo \& Elalouf~(2014)]{CovoElalouf14} Covo, S. \& Elalouf, A.~(2014). A novel single-gamma approximation to the sum of independent gamma variables, and a generalization to infinitely divisible distributions.  \emph{\eljs}\textbf{8}, 894--926.

\bibitem[Diaconis \& Zabell~(1991)]{DZ91}Diaconis, P. \& Zabell, S.~(1991). Closed Form Summation for Classical Distributions: Variations on a Theme of De Moivre.
\emph{Statist. Sci.} \textbf{6}, 284--302.

\bibitem[Eichelsbacher \& Th\"ale~(2015)]{ET15} Eichelsbacher, P. \& Th\"ale, C.~(2015). Malliavin-Stein method for Variance-Gamma approximation on Wiener space. 
\emph{\eljp}\textbf{20}, 1--28.

\bibitem[Gaunt~(2019)]{Gaunt19} Gaunt, R. E. (2019). New error bounds for Laplace approximation via Stein's method. arXiv:1911.03574.

\bibitem[Gaunt, Pickett \& Reinert~(2017)]{GPR17} Gaunt, R. E., Pickett, A. M. \& Reinert, G.~(2017). Chi-square approximation by Stein's method with application to Pearson's statistic. \emph{\anap}\textbf{27}, 720--756.

\bibitem[Goldstein \& Reinert~(1997)]{GoldsteinReinert97} Goldstein, L. \& Reinert, G.~(1997). Stein's method and the zero bias transformation with application to simple random sampling. \emph{\anap}\textbf{7}, 935--952.

\bibitem[Goldstein \& Reinert~(2005)]{GoldsteinReinert05} Goldstein, L. \& Reinert, G.~(2005). Distributional transformations, orthogonal polynomials,
and Stein characterizations. \emph{J. Theoret. Probab.}~\textbf{18}, 237--260.

\bibitem[Goldstein \& Rinott~(2005)]{GoldsteinRinott05} Goldstein, L. \& Rinott, Y.~(2005). Multivariate normal approximations by
Stein's method and size bias couplings. \emph{\jap}\textbf{33}, 1--17.

\bibitem[Goldstein \& Xia~(2006)]{Goldstein06} Goldstein, L. and Xia, A.~(2006). Zero biasing and a discrete central limit theorem. \emph{\ap}\textbf{34}, 1782--1806.

\bibitem[Kalashnikov~(1997)]{Kala97} Kalashnikov, V. V.~(1997). Geometric Sums: Bounds for Rare Events with Applications:
Risk Analysis, Reliability, Queueing. \emph{Mathematics and Its Applications}~\textbf{413}. Kluwer, Dordrecht.

\bibitem[Ledoux, Nourdin \& Peccati~(2015)]{LNP15} Ledoux, M., Nourdin, I. \& Peccati, G.~(2015). Stein's method, logarithmic Sobolev and transport inequalities. \emph{Geometric and Functional Analysis} \textbf{25}, 256--306.

\bibitem[Luk~(1994)]{Luk94} Luk, H. M.~(1994). Stein's method for the gamma distribution and related statistical applications (Doctoral dissertation, University of Southern California).

\bibitem[Nourdin \& Peccati~(2009)]{NourdinPeccati09} Nourdin, I. \& Peccati, G.~(2009). Stein's method on Wiener chaos. \emph{\ptrf}\textbf{145}, 75--118.

\bibitem[Pek\"oz \& R\"ollin~(2011)]{PekozRollin11} Pek\"oz, E. A. \& R\"ollin, A.~(2011). New rates for exponential approximation and the theorems of R\'enyi and Yaglom. \emph{\ap}\textbf{39}, 587--608.

\bibitem[R\'enyi~(1956)]{Renyi56} R\'{e}nyi, A.~(1956). A characterization of Poisson processes. Magyar Tud. Akad. Mat. Kutat\'{o}. Int. K\"{o}zl. \textbf{1}, 519--527.

\bibitem[Slepov~(2021)]{Slepov21} Slepov, N. A.~(2021). Convergence Rate of Random Geometric Sum Distributions to the Laplace Law. \emph{\tpa}\textbf{66}, 121--141.

\bibitem[Stein~(1972)]{Stein72} Stein, C.~(1972). A bound for the error in the normal approximation to the distribution of a sum of dependent random variables. \emph{In Proceedings of the Sixth Berkeley Symposium on Mathematical Statistics and Probability}, Volume 2: Probability Theory. The Regents of the University of California.

\bibitem[Sugakova~(1995)]{Sugakova95} Sugakova, E. V.~(1995). Estimates in the R\'enyi theorem for differently distributed terms. \emph{Ukrainian Mathematical Journal} \textbf{47}, 1128--1134.

\bibitem[Vallender~(1972)]{Val72} Vallender, S. S.~(1972). Calculation of the Wasserstein distance between probability distributions on the line. \emph{Theory Probab. Appl.}~\textbf{18}, 784--786. 

\bibitem[Vellaisamy \& Upadhye~(2009)] {VU09} Vellaisamy, P. \& Upadhye, N. S.~(2009). On the sums of compound negative binomial and gamma random variables. \emph{\jap}\textbf{46}, 272--283.

\end{thebibliography}
\end{document}